\documentclass[11pt]{article}
\usepackage{amssymb,epsfig}
\usepackage{amsmath}
\usepackage{graphicx}
\usepackage{epstopdf}

\title {\bf A nonlinear interpolatory reconstruction operator on non uniform grids}

\author{J.C.Trillo\thanks{
 Departamento de Matem\'atica Aplicada y Estad\'{\i}stica.
 Universidad Polit\'ecnica de Cartagena (Spain).
e-mail:{\tt jc.trillo@upct.es}} \and P. Ortiz\thanks{
 Departamento de Matem\'atica Aplicada y Estad\'{\i}stica.
 Universidad Polit\'ecnica de Cartagena (Spain).
 e-mail:{\tt poh0@alu.upct.es}} }

\renewcommand{\arraystretch}{1.8}

\begin{document}

\maketitle

\newtheorem{proposition}{Proposition}
\newtheorem{lemma}{Lemma}
\newtheorem{definition}{Definition}
\newtheorem{theorem}{Theorem}
\newtheorem{corollary}{Corollary}
\newtheorem{remark}{Remark}

\begin{abstract}
This paper is devoted to introduce the non linear reconstruction operator PPH on non uniform grids.
We define this operator and we study its main properties such as reproduction of polynomials of second degree,
approximation order and conditions for convexity preservation.
\end{abstract}

{\bf Key Words.} Interpolation,  reconstruction operator, nonlinearity,
 non uniform grids.

\vspace{10pt} {\bf AMS(MOS) subject classifications.} 41A05, 41A10, 65D17.

\section{Introduction}\label{sec1}

Reconstruction operators are widely used in computer aided geometric design. For simplicity, functions typically used as
operators are polynomials. In order to avoid undesirable phenomena generated by high degree polynomials, reconstructions are
usually built piecewise.
Due to the bad behavior of linear operators in presence of discontinuities, it has been necessary to design non linear
operators to overcome this drawback. One of these operators was defined in \cite{ADLT} and was called PPH. This operator
essentially consists of a witty modification of the classical four points piecewise Lagrange interpolation.
For the sake of simplicity, as much in theoretical analysis as in the practical implementation and computational time,
studies usually start with data given in uniform grids. Nevertheless, some applications require dealing with data over non uniform
grids. At times, it is not trivial to adapt operators defined over uniform grids to the non uniform case.
The above mentioned PPH operator was defined over an uniform grid and some of its properties were studied in \cite{ADLT}.
These reconstruction operators are the base for the definition of associated subdivision and multiresolution schemes.
In this paper, we extend the definition of the PPH reconstruction operator to data over non uniform grids and we study
some properties of this operator.\\

The paper is organized as follows: In section \ref{sec2} we recall the non linear PPH reconstruction operator \cite{ADLT}. Section \ref{sec3} is devoted
to define the PPH reconstruction operator over non uniform grids. For this purpose, the definition of the harmonic mean used in the uniform case will be
adapted to the non uniform one, so that the new reconstruction operator has a similar structure to the original PPH and amounts to it when we restrict to uniform grids. In section \ref{sec4} some properties of the PPH for non uniform
grids, such as reproduction of polynomials of second degree, approximation order, conditions for convexity preservation, will be studied. In section \ref{sec5} we introduce a modification in the harmonic mean in order to improve the approximation order for smooth functions which have a change of convexity. We present some numerical tests in section \ref{sec6}. Finally, some conclusions are included in section \ref{sec7}.

\section{Nonlinear reconstruction operator: PPH} \label{sec2}
We use this introductory section to recall the already known PPH operator on uniform grids. For more details see \cite{ADLT}.
In the first place we write the coefficients of the piecewise centered Lagrange interpolation of third degree in terms of divided differences.

Let us consider the set of points
$f_{j-1}, f_j, f_{j+1}, f_{j+2}$ corresponding to subsequent
values at the points $x_{j-1}, x_j, x_{j+1}, x_{j+2}$ of a regular grid $X$ of step $h$. Let be
$PL_j(x)$ the Lagrange interpolatory polynomial at these points, i.e.

\begin{equation}\label{eq:cond}
	PL_j(x_m) = f_m \qquad j-1\leq m \leq j+2.     
\end{equation}
$PL_j(x)$ can be expressed as
\begin{eqnarray*}
	PL_j(x) = a_0 + a_1 (x-x_{j+\frac{1}{2}}) + a_2 (x-x_{j+\frac{1}{2}})^2 + a_3 (x-x_{j+\frac{1}{2}})^3,
\end{eqnarray*}
where $x_{j+\frac{1}{2}} = \frac{x_j+x_{j+1}}{2}.$

Condition (\ref{eq:cond}) implies

\begin{equation*}
	\left(
		\begin{array}{cccc}
			1 & \left(-\frac{3}{2} h \right) & \left(-\frac{3}{2} h \right)^2 & \left(-\frac{3}{2} h \right)^3\\
			1 & \left(-\frac{1}{2} h \right) & \left(-\frac{1}{2} h \right)^2 & \left(-\frac{1}{2} h \right)^3\\
			1 & \left(\frac{1}{2} h \right)  & \left(\frac{1}{2} h \right)^2 & \left(\frac{1}{2} h \right)^3\\
			1 & \left(\frac{3}{2} h \right)  & \left(\frac{3}{2} h \right)^2 & \left(\frac{3}{2} h \right)^3\\
		\end{array}
	\right)
	\left\{
		\begin{array}{c}
			a_0\\a_1\\a_2\\a_3
		\end{array}
	\right\}=
	\left\{
		\begin{array}{c}
			f_{j-1}\\ f_j\\ f_{j+1}\\ f_{j+2}
		\end{array}
	\right\}.
\end{equation*}

Solving the system, we obtain
\begin{equation*}
	\left\{
		\begin{array}{c}
			a_0\\a_1\\a_2\\a_3
		\end{array}
	\right\}=
	\left(
		\begin{array}{cccc}
			-\dfrac{1}{16} & \dfrac{9}{16} & \dfrac{9}{16} & -\dfrac{1}{16}\\
			\dfrac{1}{24 h} & -\dfrac{9}{8 h} & \dfrac{9}{8 h} & -\dfrac{1}{24 h}\\
			\dfrac{1}{4 h^2} & -\dfrac{1}{4 h^2} & -\dfrac{1}{4 h^2} & \dfrac{1}{4 h^2}\\
			-\dfrac{1}{6 h^3} & \dfrac{1}{2 h^3} & -\dfrac{1}{2 h^3} & \dfrac{1}{6 h^3}\\
		\end{array}
	\right)
	\left\{
		\begin{array}{c}
			f_{j-1}\\ f_j\\ f_{j+1}\\ f_{j+2}
		\end{array}
	\right\}=
	\left\{
		\begin{array}{c}
			\dfrac{- f_{j-1} + 9 f_j + 9 f_{j+1} - f_{j+2}}{16}\\
			\dfrac{ f_{j-1} -27 f_j + 27 f_{j+1} - f_{j+2}}{24 h}\\
			\dfrac{ f_{j-1} - f_j - f_{j+1} + f_{j+2}}{4 h^2}\\
			\dfrac{- f_{j-1} + 3 f_j - 3 f_{j+1} + f_{j+2}}{6 h^3}
		\end{array}
	\right\}.
\end{equation*}

Introducing the divided differences
\begin{align}\label{eq:divdiff1}
	D_j &= f [ x_{j-1}, x_{j}, x_{j+1} ] = \dfrac{f_{j-1} -2 f_j + f_{j+1}}{2h^2},\notag \\
	D_{j+1} &= f [ x_{j}, x_{j+1}, x_{j+2} ] = \dfrac{f_{j} -2 f_{j+1} + f_{j+2}}{2h^2},
\end{align}
we can reformulate the coefficients of $PL_j(x)$ and the data values $f_{j-1}, \, f_{j+2}\,$ as

\begin{equation}\label{eq:Lag4coeff}
	\left\{
		\begin{array}{c}
			a_0\\a_1\\a_2\\a_3
		\end{array}
	\right\}=
	\left\{
		\begin{array}{r}
			\frac{ f_{j} + f_{j+1}}{2}-\frac{h^2}{4} \frac{D_j+D_{j+1}}{2}\\
			\frac{ - f_{j} + f_{j+1}}{h}+\frac{h}{6} (D_j-\frac{D_j+D_{j+1}}{2})\\
			\frac{D_j+D_{j+1}}{2}\\
			-\frac{2}{3h}(D_j- \frac{D_j+D_{j+1}}{2})\\
		\end{array}
	\right\},
\end{equation}

\begin{equation*}
	\begin{array}{l}
		f_{j-1} = (f_{j} + f_{j+1} - f_{j+2}) + 4h^2 \frac{D_j+D_{j+1}}{2},\\
		f_{j+2} = (- f_{j-1} + f_{j} + f_{j+1}) + 4h^2 \frac{D_j+D_{j+1}}{2}.
	\end{array}
\end{equation*}

The prediction value $f_{j+\frac{1}{2}}$ at the mid point of $x_{j},\, x_{j+1}$ of the Lagrange interpolatory polynomial $PL_j(x)$ will be
$$f_{j+\frac{1}{2}} = a_0.$$

Any smooth function $f \in C^{4}$ which passes through the points $f_{j-1}, f_j, f_{j+1}, f_{j+2}$ satisfies
$$f(x) = PL_j(x) + O(h)^{4}.$$

The corresponding prediction operator is then said of fourth accuracy.

In this case $D_j=O(1),\; D_{j+1}= O(1),\; \frac{D_j-D_{j+1}}{2}= O(h) \,$ and $\; a_1= O(1)$.

In the presence of isolated singularities, predictions made using Lagrange reconstruction operators
lose their accuracy. In fact, if $f$ has a discontinuity point in $[x_{j+1}, x_{j+2}]$, any divided difference based on a set of $m+1$ points containing $[x_{j+1}, x_{j+2}]$ verifies

\begin{equation*}
	f[x_{0}, \ldots, x_{m}]=\frac{O([f])}{h^m},
\end{equation*}\\
where $[f]=\left| f_{j+2}-f_{j+1} \right|$. So, $\; D_j=O(1) \;$ and $\; D_{j+1}= O([f])/h^2 \;$ which leads to $\frac{D_j-D_{j+1}}{2}= O([f])/h^2 \,$ and $\, a_1=O([f])/h$.

To reduce this lost of accuracy due to the presence of a singularity in $[x_{j+1}, x_{j+2}]$, in \cite{ADLT} it is proposed to replace expression of coefficient $a_1$ in (\ref {eq:Lag4coeff}) for

\begin{equation*}
	\widetilde{a}_1=\frac{ - f_{j} + f_{j+1}}{h}+\frac{h}{6} (D_j-\widetilde{D}_j),\\
\end{equation*}\\
where

\begin{equation}\label{eq:modifiedharmonic}
	\widetilde{D}_j=\left\{
		\begin{array}{cl}
			\frac{2 D_{j}D_{j+1}}{D_{j}+D_{j+1}} & \text{if} \; D_{j}D_{j+1}>0,\\
			0 & otherwise.
		\end{array}
	\right\},
\end{equation}\\
and $\frac{2 D_{j}D_{j+1}}{D_{j}+D_{j+1}}$ is the harmonic mean of $ D_{j} \,\,$ and $\,\, D_{j+1}$.

This change leads to a new reconstruction operator called PPH (Piecewise Polynomial
Harmonic) defined in \cite{ADLT} as follows

\begin{equation}\label{eq:condPolPPH}
	\left\{
		\begin{array}{lcr}
			PPH_j(x_m) & = &f_m, \qquad  j-1\leq m \leq j+1,\\
			\widetilde{a}_1 & = & \frac{ - f_{j} + f_{j+1}}{h}+\frac{h}{6} (D_j-\widetilde{D}_j).
		\end{array}
	\right\},
\end{equation}\\
i.e. $ PPH_j(x) $ is not interpolating the point affected by the singularity.

The key point of this change is that the harmonic mean is bounded as follows

\begin{equation}\label{eq:PPHKey}
	\left|
		\dfrac{2 D_{j}D_{j+1}}{D_{j}+D_{j+1}}
	\right|
	\leq
	2 min(\left|D_j\right|, \left|D_{j+1}\right|)=O(1),
\end{equation}
and hence $\widetilde{a}_1=O(1)$, instead of $O([f]/h)$ as in the linear case.

The new polynomial will be
\begin{equation}\label{eq:pphPolynomial}
	PPH_j(x) = \widetilde{a}_0 + \widetilde{a}_1 (x-x_{j+\frac{1}{2}}) + \widetilde{a}_2 (x-x_{j+\frac{1}{2}})^2 + \widetilde{a}_3 (x-x_{j+\frac{1}{2}})^3.
\end{equation}

The coefficient $\widetilde{a}_1$ is given in (\ref{eq:condPolPPH}) and the remaining ones are obtained imposing the three interpolating conditions in (\ref{eq:condPolPPH}).

Solving the system, we obtain

\begin{equation}\label{eq:pphcoeff}
	\left\{
		\begin{array}{c}
			\widetilde{a}_0\\\widetilde{a}_1\\\widetilde{a}_2\\\widetilde{a}_3
		\end{array}
	\right\}=
	\left\{
		\begin{array}{r}
			\frac{ f_{j} + f_{j+1}}{2}-\frac{h^2}{4} \widetilde{D}_j\\
			\frac{ - f_{j} + f_{j+1}}{h}+\frac{h}{6} (D_j-\widetilde{D}_j)\\
			\widetilde{D}_j\\
			-\frac{2}{3h}(D_j- \widetilde{D}_j)\\
		\end{array}
	\right\},
\end{equation}
and
\begin{equation*}
	\widetilde{f}_{j+2}= -f_{j-1}+ f_{j}+ f_{j+1}+4h^2\widetilde{D}_j,
\end{equation*}
 which can also be expressed as

\begin{equation}\label{eq:diforden}
		\widetilde{f}_{j+2}= f_{j+2}-4h^2\left(\frac{D_j+D_{j+1}}{2}-\widetilde{D}_j\right).
\end{equation}

Comparing the coefficients of $PL_j(x)$ and $PPH_j(x)$ we can see that

\begin{equation}\label{eq:difcoef}
	\left\{
		\begin{array}{c}
			\widetilde{a}_0 - a_0 \\
			\widetilde{a}_1 - a_1 \\
			\widetilde{a}_2 - a_2 \\
			\widetilde{a}_3 - a_3 \\
		\end{array}
	\right\}=
	\left\{
		\begin{array}{l}
			\frac{h^2}{4}\left(\frac{D_j+D_{j+1}}{2}-\widetilde{D}_j\right)\\
			\frac{h}{6}\left(\frac{D_j+D_{j+1}}{2}-\widetilde{D}_j\right)\\
			-\left(\frac{D_j+D_{j+1}}{2}-\widetilde{D}_j\right)\\
			-\frac{2}{3h}\left(\frac{D_j+D_{j+1}}{2}-\widetilde{D}_j\right)\\
		\end{array}
	\right\},
\end{equation}\\
i.e. $\widetilde{a}_i - a_i = O(h^{4-i})$ for $i\in [0,3]$ in smooth convex regions.

Equations (\ref{eq:diforden}) and (\ref{eq:difcoef}) show that the net effect in the coefficients $a_i$ and $f_{j+2}$ is the replacement of the arithmetic mean of $D_j$ and $D_{j+1}$ by the modified harmonic mean $\widetilde{D}_j$ in (\ref{eq:modifiedharmonic}).

\section{A nonlinear Lagrange interpolation procedure on a non uniform grid}\label{sec3}

Let us consider the set of points
$f_{j-1}, f_j, f_{j+1}, f_{j+2}$ corresponding to subsequent
values at the points $x_{j-1}, x_j, x_{j+1}, x_{j+2}$ of a non uniform grid $X$. Let be $h_{j}=x_{j}-x_{j-1}, \; h_{j+1}=x_{j+1}-x_{j}, \; h_{j+2}=x_{j+2}-x_{j+1} \;$ and $PL_j(x)$ the Lagrange interpolatory polynomial at these points, i.e.

\begin{equation}\label{eq:condNU}
	PL_j(x_m) = f_m \quad j-1\leq m \leq j+2.     
\end{equation}
$PL_j(x)$ can be expressed as
\begin{equation}
	PL_j(x) = a_0 + a_1 (x-x_{j+\frac{1}{2}}) + a_2 (x-x_{j+\frac{1}{2}})^2 + a_3 (x-x_{j+\frac{1}{2}})^3,
\end{equation}
where $x_{j+\frac{1}{2}} = \frac{x_j+x_{j+1}}{2}.$\\

Condition (\ref{eq:condNU}) implies\\
\begin{equation*}
	\left(
		\begin{array}{cccc}
			1 & \left(-h_{j}-\frac{h_{j+1}}{2} \right) & \left(-h_{j}-\frac{h_{j+1}}{2}\right)^2 & \left(-h_{j}-\frac{h_{j+1}}{2}\right)^3 \\
			1 & -\frac{h_{j+1}}{2} & \frac{h_{j+1}^2}{4} & -\frac{h_{j+1}^3}{8} \\
			1 & \frac{h_{j+1}}{2} & \frac{h_{j+1}^2}{4} & \frac{h_{j+1}^3}{8} \\
			1 & \left(\frac{h_{j+1}}{2}+h_{j+2} \right) & \left(\frac{h_{j+1}}{2}+h_{j+2}\right)^2 & \left(\frac{h_{j+1}}{2}+h_{j+2}\right)^3 \\
		\end{array}
	\right)
	\left\{
		\begin{array}{c}
			a_0\\a_1\\a_2\\a_3
		\end{array}
	\right\}=
	\left\{
		\begin{array}{c}
			f_{j-1}\\ f_j\\ f_{j+1}\\ f_{j+2}
		\end{array}
	\right\}.
\end{equation*}

Solving the system, we obtain\\
\begin{equation}\label{eq:coeffcNU}
	\left\{
		\begin{array}{c}
			a_0\\a_1\\a_2\\a_3
		\end{array}
	\right\}=
	\left(
		\begin{array}{cccc}
			c_{11} & c_{12} & c_{13} & c_{14}\\
			c_{21} & c_{22} & c_{23} & c_{24}\\
			c_{31} & c_{32} & c_{33} & c_{34}\\
			c_{41} & c_{42} & c_{43} & c_{44}\\
		\end{array}
	\right)
	\left\{
		\begin{array}{c}
			f_{j-1}\\ f_j\\ f_{j+1}\\ f_{j+2}
		\end{array}
	\right\},
\end{equation}\\
where

\renewcommand{\arraystretch}{2.5}
$
\begin{array}{ll}
	c_{11}= -\frac{h_{j+1}^2 (h_{j+1}+2 h_{j+2})}{8 h_{j} (h_{j}+h_{j+1}) (h_{j}+h_{j+1}+h_{j+2})},&
	c_{12}= \frac{(2 h_{j}+h_{j+1}) (h_{j+1}+2 h_{j+2})}{8 h_{j} (h_{j+1}+h_{j+2})},\\
	c_{13}= \frac{(2 h_{j}+h_{j+1}) (h_{j+1}+2 h_{j+2})}{8 (h_{j}+h_{j+1}) h_{j+2}},&
	c_{14}= -\frac{h_{j+1}^2 (2 h_{j}+h_{j+1})}{8 h_{j+2} (h_{j+1}+{h_{j+2}) (h_{j}+h_{j+1}+h_{j+2})}},\\
	c_{21}=  \frac{h2^2}{4 h_{j} (h_{j}+h_{j+1}) (h_{j}+h_{j+1}+h_{j+2})},&
	c_{22}=  -\frac{h_{j+1}^2+4 h_{j} (h_{j+1}+h_{j+2})}{4 h_{j} h_{j+1} (h_{j+1}+h_{j+2})},\\
	c_{23}=  \frac{h_{j+1}^2+4 h_{j+2} h_{j+1}+4 h_{j} h_{j+2}}{4 h_{j+2} h_{j+1}^2+4 h_{j} h_{j+2} h_{j+1}},&
	c_{24}= -\frac{h_{j+1}^2}{4 h_{j+2} (h_{j+1}+h_{j+2}) (h_{j}+h_{j+1}+h_{j+2})},\\
	c_{31}=  \frac{h_{j+1}+2 h_{j+2}}{2 h_{j} (h_{j}+h_{j+1}) (h_{j}+h_{j+1}+h_{j+2})},&
	c_{32}= -\frac{-2 h_{j}+h_{j+1}+2 h_{j+2}}{2 h_{j} h_{j+1}^2+2 h_{j} h_{j+2} h_{j+1}},\\
	c_{33}= -\frac{2 h_{j}+h_{j+1}-2 h_{j+2}}{2 h_{j+2} h_{j+1}^2+2 h_{j} h_{j+2} h_{j+1}},&
	c_{34}= \frac{2 h_{j}+h_{j+1}}{2 h_{j+2} (h_{j+1}+h_{j+2}) (h_{j}+h_{j+1}+h_{j+2})},\\
	c_{41}=  -\frac{1}{h_{j} (h_{j}+h_{j+1}) (h_{j}+h_{j+1}+h_{j+2})},&
	c_{42}= \frac{1}{h_{j} h_{j+1}^2+h_{j} h_{j+2} h_{j+1}},\\
	c_{43}= -\frac{1}{h_{j+2} h_{j+1}^2+h_{j} h_{j+2} h_{j+1}},&
	c_{44}=  \frac{1}{h_{j+2} (h_{j+1}+h_{j+2}) (h_{j}+h_{j+1}+h_{j+2})}.
\end{array}
$\\
\renewcommand{\arraystretch}{1.8}
\\

Introducing the divided differences
\begin{equation*}
	\begin{aligned}
		D_{j}& = & f [ x_{j-1}, x_{j}, x_{j+1} ]& = &\frac{f_{j-1}}{h_{j} (h_{j}+h_{j+1})}-\frac{f_{j}}{h_{j} h_{j+1}}+\frac{f_{j+1}}{h_{j+1} (h_{j}+h_{j+1})},\\
		D_{j+1}& = & f [ x_{j}, x_{j+1}, x_{j+2} ]& = &\frac{f_{j}}{h_{j+1} (h_{j+1}+h_{j+2})}-\frac{f_{j+1}}{h_{j+1} h_{j+2}}+\frac{f_{j+2}}{h_{j+2} (h_{j+1}+h_{j+2})},\\
	\end{aligned}
\end{equation*}
and realizing that for the uniform case, $a_{2}\,$ (equation \ref{eq:Lag4coeff}) is the arithmetic
mean of $D_{j}$ and $D_{j+1}$, we will look for a generalized mean in the expression
of $a_{2}$ in the non uniform case. This expression is

\begin{equation}\label{eq:MjNUE1}
	a_{2}= M_{j}= w_{j}D_{j}+w_{j+1}D_{j+1},
\end{equation}
where

\begin{equation}\label{eq:wNUE}
	\begin{array}{lcrcr}
		w_{j}&=&\frac{h_{j+1}+2h_{j+2}}{2(h_{j}+h_{j+1}+h_{j+2})},&&\\
		w_{j+1}&=&\frac{h_{j+1}+2h_{j}}{2(h_{j}+h_{j+1}+h_{j+2})}&=&1-w_{j}.
	\end{array}
\end{equation}

Developing equations (\ref{eq:coeffcNU}) and searching for generalized expressions of coefficients in (\ref{eq:Lag4coeff}), we obtain the whole set of coefficients of the Lagrange interpolatory polynomial $PL_{j}(x)$ in the non uniform case

\begin{equation}\label{eq:aiLagNU1}
	\left\{
		\begin{array}{c}
			a_0\\a_1\\a_2\\a_3
		\end{array}
	\right\}=
	\left\{
		\begin{array}{r}
			\frac{ f_{j} + f_{j+1}}{2}-\frac{h_{j+1}^2}{4} M_j\\
			\frac{ - f_{j} + f_{j+1}}{h_{j+1}}+\frac{h_{j+1}^2}{2(2h_{j}+h_{j+1})} (D_j-M_j)\\
			M_j\\
			-\frac{2}{2h_{j}+h_{j+1}}(D_j- M_j)\\
		\end{array}
	\right\},
\end{equation}\\
which can also be expressed as\\
\begin{equation}\label{eq:aiLagNU2}
	\left\{
		\begin{array}{c}
			a_0\\a_1\\a_2\\a_3
		\end{array}
	\right\}=
	\left\{
		\begin{array}{r}
			\frac{ f_{j} + f_{j+1}}{2}-\frac{h_{j+1}^2}{4} M_j\\
			\frac{ - f_{j} + f_{j+1}}{h_{j+1}}+\frac{h_{j+1}^2}{2(2h_{j+2}+h_{j+1})} (-D_{j+1}+M_j)\\
			M_j\\
			-\frac{2}{2h_{j+2}+h_{j+1}}(-D_{j+1}+ M_j)\\
		\end{array}
	\right\}.
\end{equation}\\

Equation (\ref{eq:MjNUE1}) can be reformulated as
\begin{equation}\label{eq:MjNUE2}
	M_{j} = k_{j-1} f_{j-1} + k_{j} f_{j} + k_{j+1} f_{j+1} + k_{j+2} f_{j+2},
\end{equation}
where the constants $\, k_{m}, \; m\in [j-1, j+2] \,$ are easily calculated. In particular,

\begin{equation}\label{eq:CoeffK}
	\begin{array}{l}
		k_{j-1} = \frac{2h_{j+2}+h_{j+1}}{2h_{j}(h_{j+1}+h_{j})(h_j+h_{j+1}+h_{j+2})},\\
		k_{j+2} = \frac{2h_j+h_{j+1}}{2h_{j+2}(h_{j+1}+h_{j+2})(h_j+h_{j+1}+h_{j+2})}.
	\end{array}
\end{equation}

Isolating the data values $f_{j-1}$ \, and \, $f_{j+2}$ results
\begin{subequations}\label{eq:fEnds}
	\begin{align}
		f_{j-1} &= \frac{-1}{k_{j-1}}(k_{j}f_{j}+k_{j+1}f_{j+1} +
			k_{j+2}f_{j+2})+ \frac{M_{j}}{k_{j-1}}, \label{eq:fEndsL}\\
		f_{j+2} &= \frac{-1}{k_{j+2}}(k_{j-1}f_{j-1}+k_{j}f_{j} +
			k_{j+1}f_{j+1})+ \frac{M_{j}}{k_{j+2}}. \label{eq:fEndsR}
	\end{align}
\end{subequations}

The prediction value $f_{j+\frac{1}{2}}$ at the mid point of $x_{j},\, x_{j+1}$ of the Lagrange
interpolatory polynomial $PL_j(x)$ will be
$$f_{j+\frac{1}{2}} = a_0.$$

Any smooth function $f \in C^{4}$ which passes through the points $f_{j-1}, f_j, f_{j+1}, f_{j+2}$ satisfies
$$f(x) = PL_j(x) + O(h)^{4},$$
where $h = \max\{h_j, h_{j+1}, h_{j+2}\}.$\\
The corresponding prediction operator is then said of fourth accuracy. In this case
$D_j=O(1),\, D_{j+1}= O(1)\,$ and $\, D_j-M_j= w_{j+1}(D_j-D_{j+1})= O(h)$.

In the presence of isolated singularities, predictions made using Lagrange reconstruction operators
lose their accuracy as much in the uniform as in the non uniform case.

In order to avoid the undesirable effects of discontinuities, we introduce the weighted harmonic mean over non uniform grids,
which will be used in the definition of the extension of the PPH reconstruction operator to such grids.
This harmonic mean is built as the inverse of the weighted arithmetic mean of the inverses of the given values, i.e.

\begin{equation}
    \widetilde{V}(x,y)= \frac{1}{w_x \frac{1}{x}+w_y \frac{1}{y}} = \frac{x y}{w_x y + w_y x},
\end{equation}\\
where $\; w_x, \, w_y \;$ stands for the considered weights.\\
The key points of this change are:
\begin{itemize}
	\item[1.] If $x\geq 0$, $y \geq 0$, and $\; x=\min\{x, y\}\;$ the harmonic mean is bounded as follows
	\begin{equation} \label{eq:boundedVnu}
		|\widetilde{V}(x,y)| < \min\left\{ \frac{1}{w_x}|x|,\frac{1}{w_y}|y|\right\} \leq \frac{1}{w_x}|x|.
	\end{equation}

	\item[2.] If $x= O(1)$ and $y= O(1)$ are close to each other, i.e. $\, |x-y|= O(h)$, then the weighted harmonic mean is also
close to the weighted arithmetic mean,
	\begin{equation}\label{eq:distMVnu}
		|M(x, y)-\widetilde{V}(x,y)|= \frac{w_x w_y}{w_x y+w_y x}(x-y)^2= O(h^2).
	\end{equation}
\end{itemize}

Notice that all these expressions coincide with their equivalent ones in the uniform case when $h_j = h_{j+1} = h_{j+2} \;\; \forall j.$\\

When a singularity appears outside the central interval, i.e. at $[x_{j-1}, x_{j} ]$ or at $[x_{j+1}, x_{j+2}]$, we propose
not to interpolate the data at such a point and change it for another value. There
are two proposed changes, depending on the interval where the singularity lies.\\

\textbf{Case 1. $|D_j| \leq |D_{j+1}|$}, i.e, the possible singularity is at $[x_{j+1}, x_{j+2}] $. We propose to replace $f_{j+2}$
by $\widetilde{f}_{j+2}$ changing the weighted arithmetic mean in
equation (\ref{eq:fEndsR}) by the weighted harmonic mean, as follows

\begin{equation}\label{eq:fRightEnd}
	\widetilde{f}_{j+2} = \frac{-1}{k_{j+2}}(k_{j-1}f_{j-1}+k_{j}f_{j}+k_{j+1}f_{j+1})+ \frac{\widetilde{V}_j}{k_{j+2}},
\end{equation}
where
\begin{equation}\label{eq:Vwidetilde}
    \widetilde{V}_j=
        \left\{
            \begin{array}{ll}
                \widetilde{V}(D_j, D_{j+1})=\frac{D_j D_{j+1}}{w_j D_{j+1} + w_{j+1} D_j} & if D_{j} D_{j+1} > 0,\\
                0 & otherwise.
            \end{array}
        \right\}
\end{equation}

So,
\begin{equation}\label{eq:distfftildeNU}
    |\widetilde{f}_{j+2}-f_{j+2}|= \frac{2h_{j+2}(h_{j+1}+h_{j+2})(h_j+h_{j+1}+h_{j+2})}{2h_j+h_{j+1}}|M_j-\widetilde{V}_j|.
\end{equation}

Notice that in smooth regions, taking into account equation (\ref{eq:distMVnu}) we have
$\; |\widetilde{f}_{j+2}-f_{j+2}|=O(h^4)$.

It is also important to point out that equation (\ref{eq:fRightEnd}) says that
$\; \widetilde{f}_{j+2} \;$ is not significantly affected by the possible singularity
at the interval $\; [x_{j+1}, x_{j+2}]\;$ since by property
(\ref{eq:boundedVnu}) $\; |\widetilde{V}_j| \leq \frac{1}{w_j}|D_j|\;$ and in turn
$\; D_j \;$ is not affected by such discontinuity.\\

Previous change leads to a new reconstruction operator PPH which is an
extension to non uniform grids of the PPH (Piecewise Polynomial Harmonic)
and is defined by the following conditions

\begin{equation}\label{eq:condPolPPHnu}
	\left\{
		\begin{array}{lcl}
			PPH_j(x_m) & = &f_m, \qquad  j-1\leq m \leq j+1,\\
			PPH_j(x_{j+2}) & = & \widetilde{f}_{j+2}.
		\end{array}
	\right\}
\end{equation}\\

The new polynomial takes the form
\begin{equation}\label{eq:pphPolynomialNU}
	PPH_j(x) = \widetilde{a}_0 + \widetilde{a}_1 (x-x_{j+\frac{1}{2}}) + \widetilde{a}_2 (x-x_{j+\frac{1}{2}})^2 + \widetilde{a}_3 (x-x_{j+\frac{1}{2}})^3.
\end{equation}

Imposing conditions (\ref{eq:condPolPPHnu}) in (\ref{eq:pphPolynomialNU}) lead us to the following linear system

\begin{equation*}
	\left(
		\begin{array}{cccc}
			1 & \left(-h_{j}-\frac{h_{j+1}}{2} \right) & \left(-h_{j}-\frac{h_{j+1}}{2}\right)^2 & \left(-h_{j}-\frac{h_{j+1}}{2}\right)^3 \\
			1 & -\frac{h_{j+1}}{2} & \frac{h_{j+1}^2}{4} & -\frac{h_{j+1}^3}{8} \\
			1 & \frac{h_{j+1}}{2} & \frac{h_{j+1}^2}{4} & \frac{h_{j+1}^3}{8} \\
			1 & \left(\frac{h_{j+1}}{2}+h_{j+2} \right) & \left(\frac{h_{j+1}}{2}+h_{j+2}\right)^2 & \left(\frac{h_{j+1}}{2}+h_{j+2}\right)^3 \\
		\end{array}
	\right)
	\left\{
		\begin{array}{c}
			\widetilde{a}_0\\\widetilde{a}_1\\\widetilde{a}_2\\\widetilde{a}_3
		\end{array}
	\right\}=
	\left\{
		\begin{array}{c}
			f_{j-1}\\ f_j\\ f_{j+1}\\ \widetilde{f}_{j+2}
		\end{array}
	\right\}.
\end{equation*}

Solving it, we obtain

\begin{equation}\label{eq:pphcoeffNUr}
	\left\{
		\begin{array}{c}
			\widetilde{a}_0\\\widetilde{a}_1\\\widetilde{a}_2\\\widetilde{a}_3
		\end{array}
	\right\}=
	\left\{
		\begin{array}{r}
			\frac{ f_{j} + f_{j+1}}{2}-\frac{h_{j+1}^2}{4} \widetilde{V}_j\\
			\frac{ - f_{j} + f_{j+1}}{h_{j+1}}+\frac{h_{j+1}^2}{4 h_{j}+ 2 h_{j+1}} (D_j-\widetilde{V}_j)\\
			\widetilde{V}_j\\
			-\frac{2}{2h_{j}+h_{j+1}}(D_j- \widetilde{V}_j)\\
		\end{array}
	\right\}.
\end{equation}\\

The difference of these coefficients with the ones of $PL_j(x)$ given in equation (\ref{eq:aiLagNU1}) is

\begin{equation}\label{eq:difcoefNU1}
	\left\{
		\begin{array}{c}
			\widetilde{a}_0 - a_0 \\
			\widetilde{a}_1 - a_1 \\
			\widetilde{a}_2 - a_2 \\
			\widetilde{a}_3 - a_3 \\
		\end{array}
	\right\}=
	\left\{
		\begin{array}{l}
			\frac{h_{j+1}^2}{4}\left(M_j-\widetilde{V}_{j}\right)\\
			\frac{h_{j+1}^2}{4 h_{j}+ 2 h_{j+1}} (M_j-\widetilde{V}_j)\\
			-(M_j-\widetilde{V}_j)\\
			-\frac{2}{2h_{j}+h_{j+1}}(M_j- \widetilde{V}_j)\\
		\end{array}
	\right\}.
\end{equation}

These expressions will be used in next section to prove some properties of the defined operators.\\

\textbf{Case 2. $|D_j| > |D_{j+1}|$}, i.e, the possible singularity is at  $[x_{j-1}, x_{j}] $.
We propose to replace $f_{j-1}$ by $\widetilde{f}_{j-1}$ changing the weighted arithmetic mean in equation (\ref{eq:fEndsL}) by the weighted harmonic mean.\\

By developing this case in a similar way to the previous one, we obtain the following
coefficients for the polynomial (\ref{eq:pphPolynomialNU})

\begin{equation}\label{eq:pphcoeffNUl}
	\left\{
		\begin{array}{c}
			\widetilde{a}_0\\\widetilde{a}_1\\\widetilde{a}_2\\\widetilde{a}_3
		\end{array}
	\right\}=
	\left\{
		\begin{array}{r}
			\frac{ f_{j} + f_{j+1}}{2}-\frac{h_{j+1}^2}{4} \widetilde{V}_j\\
			\frac{ - f_{j} + f_{j+1}}{h_{j+1}}+\frac{h_{j+1}^2}{2 h_{j+1}+4 h_{j+2}} (-D_{j+1}+\widetilde{V}_j)\\
			\widetilde{V}_j\\
			-\frac{2}{h_{j+1}+2h_{j+2}}(-D_{j+1}+ \widetilde{V}_j)\\
		\end{array}
	\right\}.
\end{equation}\\
 Their difference with the coefficients of $PL_j(x)$ given in equation (\ref{eq:aiLagNU2}) is in this case\\
  \begin{equation}\label{eq:difcoefNU2}
	\left\{
		\begin{array}{c}
			\widetilde{a}_0 - a_0 \\
			\widetilde{a}_1 - a_1 \\
			\widetilde{a}_2 - a_2 \\
			\widetilde{a}_3 - a_3 \\
		\end{array}
	\right\}=
	\left\{
		\begin{array}{l}
			\frac{h_{j+1}^2}{4}\left(M_j-\widetilde{V}_{j}\right)\\
			-\frac{h_{j+1}^2}{2 h_{j+1}+4 h_{j+2}} (M_j-\widetilde{V}_j)\\
			-(M_j-\widetilde{V}_j)\\
			\frac{2}{2h_{j+2}+h_{j+1}}(M_j- \widetilde{V}_j)\\
		\end{array}
	\right\}.
\end{equation}\\

\section{Properties}\label{sec4}
In this section we study some interesting properties of the new reconstruction operator. We start with the property of reproduction of polynomials up to degree two.
\subsection{Reproduction of polynomials of degree 2}
If our function $f(x)$ is a polynomial of degree 2, then $D_{j}=D_{j+1}=D$ and $D_{j}D_{j+1}=D^2>0$.
 Using equations (\ref{eq:MjNUE1}), (\ref{eq:Vwidetilde}) and (\ref{eq:difcoefNU1}) we get
\begin{align*}
	M_{j} &= w_{j}D+(1-w_{j})D = D, \\
	\widetilde{V}_{j} &=\frac{D^2}{w_j D + (1-w_{j}) D} = D,\\
    \widetilde{a}_i &= a_i \;\; \forall i\in [0,3].
 \end{align*}
 So, $PPH_j(x) = PL_j(x)$, i.e. $PPH_j(x)$ reproduce polynomials of degree 2, since $PL_j(x)$ does it.

\subsection{Approximation Order}\label{sec42}
From equations (\ref{eq:MjNUE1}) and (\ref{eq:Vwidetilde}) we can write
\begin{equation}\label{eq:difMV}
    M_{j}-\widetilde{V}_j=
        \left\{
            \begin{array}{ll}
                \dfrac{w_{j} w_{j+1} (D_{j+1} - D_{j})^2}{w_j D_{j+1} + w_{j+1} D_j} & if D_{j} D_{j+1} > 0,\\
                M_{j} & otherwise.
            \end{array}
        \right\}
\end{equation}

We consider two cases. The first one when the function $\,f\,$ is smooth and the second one when $\,f\,$ has a singularity
outside the central interval\\
\\
\textbf{Case 1}. Function $f \in C^4$ is smooth, i.e.
$$D_{j}= O(1), \; D_{j+1}= O(1) \; \mbox{and} \; D_{j+1} - D_{j}= O(h),$$
where $h = \max\{h_j, h_{j+1}, h_{j+2}\}.$\\
\\
\textbf{Case 1.1}. $D_{j} D_{j+1} > 0$.
\\
\\
Equations (\ref{eq:difcoefNU1}), (\ref{eq:difcoefNU2}) and (\ref{eq:difMV}) let us write
\begin{align*}
	|M_{j} - \widetilde{V}_{j}| &= O(h^2), \\
    |\widetilde{a}_i - a_i| &=  O(h^{4-i}) \quad \forall i \in [0, 3],
\end{align*}
so,
\begin{equation*}
|PPH_{j}(x)-PL_{j}(x)| \leq \sum_{i=0}^{3} |\widetilde{a}_i - a_i||(x-x_{j+1/2})^i|= O(h^4).
\end{equation*}

Taking into account the triangular inequality
\begin{equation*}
	|f(x)-PPH_j(x)| \leq |f(x)-PL_j(x)|+|PL_j(x)-PPH_j(x)| = O(h^4),
\end{equation*}
that is,  when the function $\, f \,$ is smooth and convex in $\, [x_{j-1}, x_{j+2}] \,$ we have fourth order accuracy.
\\
\\
\textbf{Case 1.2.}  $D_{j}D_{j+1} \leq 0$.\\
\\
Using again equations (\ref{eq:difcoefNU1}), (\ref{eq:difcoefNU2}), (\ref{eq:difMV})
and the triangular inequality, we obtain
\begin{align*}
	|M_{j} - \widetilde{V}_{j}| &= O(1), \\
    |\widetilde{a}_i - a_i| &=  O(h^{2-i}) \quad \forall i \in [0, 3],\\
    |PPH_{j}(x)-PL_{j}(x)| &\leq \sum_{i=0}^{3} |\widetilde{a}_i - a_i||(x-x_{j+1/2})^i|= O(h^2),\\
    |f(x)-PPH_j(x)| &\leq |f(x)-PL_j(x)|+|PL_j(x)-PPH_j(x)| = O(h^2), \notag
\end{align*}
therefore in this case the accuracy is reduced to second order.
\\
\\
\textbf{Case 2}. Function $f$ has a singularity outside the central interval.\\
\\
Let us suppose that the singularity is at $[x_{j+1},x_{j+2}]\;$ and the function
is smooth in $[x_{j-1},x_{j+1}] \;$, i.e. $D_j=O(1)\;$ and $|D_j| \leq|D_{j+1}| \;$.\\
If the singularity is at $[x_{j-1},x_{j}]\;$ and the function is smooth in $[x_{j},x_{j+2}] \;$,
then $D_{j+1}=O(1), \; |D_j| >|D_{j+1}|\;$ and the process is similar to the one
shown below.\\

Let be $PL2_j(x)$ the second degree Lagrange interpolatory polynomial at points
$(x_{j-1}, f_{j-1}), \:(x_{j}, f_{j}), \: (x_{j+1}, f_{j+1})$.
\begin{eqnarray*}
	PL2_j(x) =\widehat{a}_0 + \widehat{a}_1 (x-x_{j+\frac{1}{2}})+ \widehat{a}_2 (x-x_{j+\frac{1}{2}})^2,
\end{eqnarray*}
where
\begin{equation}\label{eq:Lag3coeff}
	\left\{
		\begin{array}{c}
			\widehat{a}_0\\\widehat{a}_1\\\widehat{a}_2
		\end{array}
	\right\}=
	\left\{
		\begin{array}{r}
			\frac{ f_{j} + f_{j+1}}{2}-\frac{h_{j+1}^2}{4}D_j\\
			\frac{ - f_{j} + f_{j+1}}{h_{j+1}}\\
            D_j\\
		\end{array}
	\right\}.
\end{equation}

The difference between these coefficients and the ones of $PPH_j(x)$ shown in equation (\ref{eq:pphcoeffNUr}) is given by

\begin{equation}\label{eq:difcoefLag3_1}
	\left\{
		\begin{array}{c}
			\widetilde{a}_0 - \widehat{a}_0 \\
			\widetilde{a}_1 - \widehat{a}_1 \\
			\widetilde{a}_2 - \widehat{a}_2 \\
			\widetilde{a}_3\\
		\end{array}
	\right\}=
	\left\{
		\begin{array}{l}
			\frac{h_{j+1}^2}{4}\left(D_j-\widetilde{V}_{j}\right)\\
			\frac{h_{j+1}^2}{4 h_{j}+ 2 h_{j+1}} (D_j-\widetilde{V}_j)\\
			-(D_j-\widetilde{V}_j)\\
			-\frac{2}{2h_{j}+h_{j+1}}(D_j- \widetilde{V}_j)\\
		\end{array}
	\right\}.
\end{equation}\\
\\
\textbf{Case 2.1}. $D_{j} D_{j+1} > 0$.\\
\\
Taking into account equations (\ref{eq:boundedVnu}), (\ref{eq:Vwidetilde}), (\ref{eq:difcoefLag3_1}) and the triangular inequality
\begin{align*}
    |\widetilde{V}(D_j, D_{j+1})| &\leq \frac{1}{w_j}|D_j|,\\
    |D_j-\widetilde{V}_j| &\leq|D_j|+\frac{1}{w_j}|D_j|=\frac{1+w_j}{w_j}|D_j|=O(1),\\
    |\widetilde{a}_i - \widehat{a}_i| &=O(h^{2-i}) \quad \forall i \in [0, 3],\\
    |PPH_{j}(x)-PL2_{j}(x)| &\leq \sum_{i=0}^{3} |\widetilde{a}_i - \widehat{a}_i||(x-x_{j+1/2})^i|= O(h^2),\\
	|f(x)-PPH_j(x)| &\leq |f(x)-PL2_j(x)|+|PL2_j(x)-PPH_j(x)| = O(h^2).
\end{align*}
\\
\textbf{Case 2.2.}  $D_{j}D_{j+1} \leq 0$.\\
\\
Equations (\ref{eq:Vwidetilde}), (\ref{eq:difcoefLag3_1}) and the triangular inequality  lead us to
\begin{align*}
    \widetilde{V}_j &=0,\\
    |D_j-\widetilde{V}_j| &=O(1),\\
    |\widetilde{a}_i - \widehat{a}_i| &=O(h^{2-i}) \quad \forall i \in [0, 3],\\
    |PPH_{j}(x)-PL2_{j}(x)| &\leq \sum_{i=0}^{3} |\widetilde{a}_i - \widehat{a}_i||(x-x_{j+1/2})^i|= O(h^2),\\
	|f(x)-PPH_j(x)| &\leq |f(x)-PL2_j(x)|+|PL2_j(x)-PPH_j(x)| = O(h^2). \notag
\end{align*}

We observe that when the singularity is at $\, [x_{j-1}, x_{j}] \,$ or at $\, [x_{j+1}, x_{j+2},] \,$ we do not lose all accuracy, but we maintain at least second order accuracy. Unfortunately, when the singularity lies on the central interval $\, [x_{j}, x_{j+1}] \,$ this approach does not allow us to obtain any gain with respect to other reconstruction operators.

\subsection{Convexity preservation}

Let  $\, (x_{j+s}, f_{j+s}), \: -1 \leq s \leq 2 \;$ be a convex set  of points, i.e. $D_{j} D_{j+1} > 0$. Let us also consider
that $\,D_j>0, \, D_{j+1} >0.$  If $\,D_j < 0, \, D_{j+1} < 0, \,$ we can proceed in a completely similar way to the one that we will study next.

Under previous hypothesis, it is said that a reconstruction operator $\, R(x) \,$ strictly preserves convexity in the interval
$\, [a, b] \,$ if
\begin{equation}\label{eq:NUEconvexCondit1}
	R^{''}(x) > 0 \quad \forall x \in [a, b].
\end{equation}

 For the particular case of $\, PPH_{j}(x), \,$ computing derivatives in equation (\ref{eq:pphPolynomialNU}), last condition can be expressed as follows, .
\begin{equation}\label{eq:NUEderivadas1}
	PPH_{j}^{''}(x)=2 \widetilde{a}_2 + 6 \widetilde{a}_3 (x-x_{j+\frac{1}{2}}) > 0.
\end{equation}

In order to analyze the sign of $\, PPH_{j}^{''}(x) \,$ we need to consider two cases due to the fact that the expression of
$\, PPH_{j}(x) \,$ is different for  $\,|D_{j}| \leq |D_{j+1}|\,$ than for  $\,|D_{j}| > |D_{j+1}|$.
\\
\\
\textbf{Case 1}. $|D_{j}| \leq |D_{j+1}|$.
\\
\\
Replacing coefficients $\widetilde{a}_2, \widetilde{a}_3 \,$ coming from equation (\ref{eq:pphcoeffNUr}) in condition (\ref{eq:NUEderivadas1}) results
\begin{equation}\label{eq:eq2}
	PPH_{j}^{''}(x)=2 \widetilde{V}_j - \frac{12}{2h_j+h_{j+1}}(D_j-\widetilde{V}_j) (x-x_{j+\frac{1}{2}})>0.
\end{equation}

Rearranging equation (\ref{eq:Vwidetilde}) we can see that for the present case
\begin{equation*}
	D_j = w_j\widetilde{V}_j+w_{j+1}\frac{D_j}{D_{j+1}}\widetilde{V}_j \leq (w_j+w_{j+1})\widetilde{V}_j=\widetilde{V}_j,
\end{equation*}
i.e. $\; \widetilde{V}_j-D_j\geq 0$. So, the expression for the abscissas satisfying condition (\ref{eq:eq2}) is given by
\begin{equation}\label{eq:NUEconvexCondit2}
	x > x_{j+\frac{1}{2}} - \dfrac{2 h_{j}+h_{j+1}}{6} \dfrac{\widetilde{V}_j}{\widetilde{V}_j-D_j}.
\end{equation}

Replacing $\widetilde{V}_j$ for its expression in equation (\ref{eq:Vwidetilde}) and later $w_{j}, w_{j+1}$ by theirs in
(\ref{eq:wNUE}) results

\begin{equation}\label{eq:NUEconvexCondit3}
	x > x_{j+\frac{1}{2}} - \dfrac{ h_{j}+h_{j+1}+h_{j+2}}{3} \dfrac{D_{j+1}}{D_{j+1}-D_{j}}.
\end{equation}

Evaluating previous expression at $x_{j} \,$ we obtain the condition for convexity preservation in $[x_j, x_{j+2}]$. This condition reads

\begin{equation}\label{eq:NUEconvexCondit4_1}
    \left({h_{j+1}-2(h_{j}+h_{j+2})}\right)D_{j+1}  < 3 h_{j+1}D_{j}.
\end{equation}

Thus, if $\, h_{j+1} \leq 2(h_{j}+h_{j+2}) \,$ the condition is satisfied. If on the contrary $\, h_{j+1} > 2(h_{j}+h_{j+2}) \,$ we have

\begin{equation}\label{eq:NUEconvexCondit4_2}
    \dfrac{D_{j+1}}{D_{j}} < \dfrac{3 h_{j+1}}{h_{j+1}-2(h_{j}+h_{j+2})}.
\end{equation}

Evaluating expression (\ref{eq:NUEconvexCondit3}) at $x_{j-1} \,$ we obtain the condition for convexity preservation in $[x_{j-1}, x_{j+2}]$. This condition is

\begin{equation}\label{eq:NUEconvexCondit5_1}
    \left(4h_{j}+h_{j+1}-2h_{j+2}\right)D_{j+1}  < 3 (2 h_{j}+h_{j+1})D_{j}.
\end{equation}

Thus, if $\, 4h_{j}+h_{j+1} \leq 2h_{j+2} \,$ the condition is fulfilled. If on the contrary $\, 4h_{j}+h_{j+1} > 2h_{j+2} \,$ we have

\begin{equation}\label{eq:NUEconvexCondit5_2}
    \dfrac{D_{j+1}}{D_{j}} < \dfrac{3(2 h_{j}+h_{j+1})}{4 h_{j}+h_{j+1}-2h_{j+2}}.
\end{equation}

Working in a similar way with the Lagrange reconstruction operator $PL_{j}(x)$, we obtain the following analogue expressions to (\ref{eq:NUEconvexCondit2}) and (\ref{eq:NUEconvexCondit3}) for the abscissas fulfilling condition  $PL_{j}^{''}(x) > 0$

\begin{equation}\label{eq:NUEconvexConditPL2}
	x > x_{j+\frac{1}{2}} - \dfrac{2 h_{j}+h_{j+1}}{6} \dfrac{M_j}{M_j-D_j},
\end{equation}

\begin{equation}\label{eq:NUEconvexConditPL3}
	x > x_{j+\frac{1}{2}}  - \dfrac{2 h_{j}+h_{j+1}}{6} - \dfrac{ h_{j}+h_{j+1}+h_{j+2}}{3} \dfrac{D_{j}}{D_{j+1}-D_{j}}.
\end{equation}

If we call $X_{PPH}$ and $X_{PL}$  to the second member of inequalities (\ref{eq:NUEconvexCondit3}) and (\ref{eq:NUEconvexConditPL3}) respectively, their difference will be

\begin{equation}\label{eq:NUEconvexDif}
	X_{PL}-X_{PPH} =  \dfrac{h_{j+1}+ 2 h_{j+2}}{6}>0,
\end{equation}
i.e, $PPH$ reconstruction operator preserves the convexity in a wider interval than Lagrange reconstruction operator does.
\\
\\
\textbf{Case 2}. $|D_{j}| > |D_{j+1}|$.
\\
\\
Replacing coefficients $\widetilde{a}_2, \widetilde{a}_3 \,$ coming from equation (\ref{eq:pphcoeffNUl}) in equations (\ref{eq:NUEderivadas1})  and following a similar track to the shown, we obtain next expressions for the abscissas verifying $PPH_{j}^{''}(x) > 0$
\begin{equation}\label{eq:bNUEconvexCondit2}
	x < x_{j+\frac{1}{2}} + \dfrac{h_{j+1}+2 h_{j+2}}{6} \dfrac{\widetilde{V}_j}{\widetilde{V}_j-D_{j+1}},
\end{equation}

\begin{equation}\label{eq:bNUEconvexCondit3}
	x < x_{j+\frac{1}{2}} + \dfrac{ h_{j}+h_{j+1}+h_{j+2}}{3} \dfrac{D_{j}}{D_{j}-D_{j+1}},
\end{equation}
\\

and these others for the abscissas satisfying $PL_{j}^{''}(x) > 0$
\begin{equation}\label{eq:bNUEconvexConditPL2}
	x < x_{j+\frac{1}{2}} + \dfrac{h_{j+1}+ 2 h_{j+2}}{6} \dfrac{M_j}{M_j-D_{j+1}},
\end{equation}

\begin{equation}\label{eq:bNUEconvexConditPL3}
	x < x_{j+\frac{1}{2}} + \dfrac{h_{j+1}+ 2 h_{j+2}}{6} + \dfrac{ h_{j+1}+  h_{j+2}+h_{j+2}}{3} \dfrac{D_{j+1}}{D_{j}-D_{j+1}}.
\end{equation}

If we call $X_{PPH}$ and $X_{PL}$  to the second member of inequalities  (\ref{eq:bNUEconvexCondit3}) and (\ref{eq:bNUEconvexConditPL3}) respectively, we observe that

\begin{equation}\label{eq:bNUEconvexDif}
	X_{PPH}-X_{PL} =  \dfrac{2 h_{j}+ h_{j+1}}{6}>0,
\end{equation}
i.e, $PPH$ reconstruction operator preserves the convexity in a wider interval than Lagrange reconstruction operator does also in this case.

\section{Translation}\label{sec5}
As we have seen in section \ref{sec4}, if a smooth function verifies $\; D_{j} D_{j+1} \leq 0, \; PPH \;$ reconstruction gives an approximation of order $\; O(h^2) \;$ (case 1.2) lower than $\; O(h^4)$ obtained with the weighted arithmetic mean in the Lagrange case. Also even if $\; D_{j} D_{j+1} > 0, \;$ in some cases, such as near a inflection point, it could happen that for some $h, \;$ $D_j = O(1)$ and $D_{j+1} = O(1), \;$ but they are of the same magnitude or even smaller than $h$, then equation (\ref{eq:difMV}) give us $M_j-\widetilde{V}_j = O(h)$ and $f(x)-PPH_j(x) = O(h^3)$.
\\

To avoid these problems, we propose to define a new mean $J_j$ which makes use of a translation $T$ that verify next conditions:
\\
\\
1.  $(D_j+T)(D_{j+1}+T) > 0. \;$ So, if we consider the differences $\,D_j+T\,$ and $\,D_{j+1}+T\,$ instead of $D_j$ and $D_{j+1},$ case 1.2 in section \ref{sec42} is turned on case 1.1 when we work with smooth functions.
\\
\\
2.  $M_j(D_j+T, D_{j+1}+T)-\widetilde{V}_j(D_j + T, D_{j+1}+T)= O(h^2)\;$ for smooth functions. Taking into account equations (\ref{eq:MjNUE1}) and (\ref{eq:Vwidetilde}), we can see that
\begin{equation*}
	M_j(D_j+T, D_{j+1}+T)-\widetilde{V}_j(D_j + T, D_{j+1}+T)= \dfrac{w_{j} w_{j+1} (D_{j+1} -  D_{j})^2}{w_{j} D_{j+1} + w_{j+1} D_{j}+T}.
\end{equation*}

One possible definition of translation $T$ fulfilling previous conditions is
\begin{equation*}
    T=
        \left\{
            \begin{array}{ll}
                sign \left[ \max \left(|D_{j}|, |D_{j+1}| \right) \right]\varepsilon & if D_{j} D_{j+1} > 0,\\
                sign \left[ \max \left(|D_{j}|, |D_{j+1}| \right) \right]\left[ \min \left(|D_{j}|, |D_{j+1}| \right)+ \varepsilon \right] & otherwise.
            \end{array}
        \right\}
\end{equation*}
where $\, \varepsilon = O(1) \,$ is a constant.
\\

The proposed new mean is
\begin{equation}
	J_j(D_j, D_{j+1})=\widetilde{V}_j(D_j + T, D_{j+1}+T) -T,
\end{equation}
which verifies
\begin{equation*}
	\left|M_j(D_j, D_{j+1})-J_j(D_j , D_{j+1})\right|= \left|M_j(D_j+T, D_{j+1}+T)-\widetilde{V}_j(D_j + T, D_{j+1}+T)\right|=O(h^2),
\end{equation*}
even when $D_j$ and $D_{j+1}$ are of the same magnitude or smaller than $h$. From this property results
\begin{align*}
    |\widetilde{a}_i - a_i| &=O(h^{4-i}) \quad \forall i \in [0, 3],\\
    |PPH_{j}(x)-PL_{j}(x)| &\leq \sum_{i=0}^{3} |\widetilde{a}_i - a_i||(x-x_{j+1/2})^i|= O(h^4),\\
    |f(x)-PPH_j(x)| &\leq |f(x)-PL_j(x)|+|PL_j(x)-PPH_j(x)| = O(h^4). \notag
\end{align*}

The new mean also verifies this other property

\begin{equation}\label{eq:bNUETransP2}
	|J_j(D_j, D_{j+1})| \leq
        \left\{
            \begin{array}{ll}
                \dfrac{|D_{j}|}{w_{j}} + \dfrac{w_{j+1}}{w_{j}}\varepsilon & if |D_{j}| \leq |D_{j+1}|,\\
                \dfrac{|D_{j+1}|}{w_{j+1}} + \dfrac{w_{j}}{w_{j+1}}\varepsilon & otherwise.
            \end{array}
        \right\}
\end{equation}
which allows us to prove the results of previous case 2 about approximation order in the presence of a singularity (section \ref{sec42}) following the same track, that is
\begin{equation*}
	|f(x)-PPH_j(x)| \leq |f(x)-PL2_j(x)|+|PL2_j(x)-PPH_j(x)| = O(h^2).
\end{equation*}


The election of $\varepsilon, \;$ is a compromise between two previous properties, so that the increase of $\varepsilon \;$ closes the reconstruction to the function in smooth areas but reduces the approach near singularities.\\

In order to demonstrate property (\ref{eq:bNUETransP2}), let us suppose that the function
is smooth in $[x_{j-1},x_{j+1}] \;$, i.e. $D_j=O(1)\;$ and $|D_j| \leq|D_{j+1}|\;$. The singularity, if any, will be at $[x_{j+1},x_{j+2}]$.

If the function is smooth in $[x_{j},x_{j+2}] \;$, i.e. $D_{j+1}=O(1),\;$ $|D_j| > |D_{j+1}|\;$ and the possible singularity is at $[x_{j-1},x_{j}]\;$ we follow a similar process to the one shown below.\\

With previous supposition there are four possible cases\\
\\
\textbf{Case A. $D_j \leq 0, D_{j+1}>0$}. In this case $T=- D_{j}+\varepsilon > 0.$

\begin{equation*}
	J_j(D_j, D_{j+1})=\widetilde{V}_j(\varepsilon, D_{j+1}- D_{j}+\varepsilon) +  D_{j}-\varepsilon,
\end{equation*}

We know that
\begin{align*}
	\widetilde{V}_j(\varepsilon, D_{j+1}- D_{j}+\varepsilon) &\geq \varepsilon,\\
	\widetilde{V}_j(\varepsilon, D_{j+1}- D_{j}+\varepsilon) &< \dfrac{1}{w_j}(\varepsilon).
\end{align*}

Thus, if 	$\; \widetilde{V}_j(\varepsilon, D_{j+1}- D_{j}+\varepsilon) + D_{j}-\varepsilon \geq 0, $

\begin{equation*}
	|J_{j}(D_j, D_{j+1})| < \left|D_{j}\right|+ \dfrac{w_{j+1}}{w_{j}} \varepsilon = O(1).
\end{equation*}

If  $\; \widetilde{V}_j(\varepsilon, D_{j+1}- D_{j}+\varepsilon) + D_{j}-\varepsilon < 0, $

\begin{equation*}
	|J_{j}(D_j, D_{j+1})| = \varepsilon - D_{j}  - \widetilde{V}_j(\varepsilon, D_{j+1}- D_{j}+\varepsilon) \leq |D_{j}| = O(1).
\end{equation*}
\\
\textbf{Case B. $D_j \geq 0, \, D_{j+1}<0$}. In this case $T=-D_{j}-\varepsilon < 0.$

\begin{equation*}
	J_{j}(D_j, D_{j+1})=\widetilde{V}_j(- \varepsilon, D_{j+1}- D_{j}-\varepsilon) +  D_{j}+\varepsilon,
\end{equation*}

We know that
\begin{align*}
	\widetilde{V}_j(- \varepsilon, D_{j+1}- D_{j}-\varepsilon) &\leq - \varepsilon,\\
	\widetilde{V}_j(- \varepsilon, D_{j+1}- D_{j}-\varepsilon) &> \dfrac{1}{w_j}( -\varepsilon).
\end{align*}

Thus, if $\; \widetilde{V}_j(- \varepsilon, D_{j+1}- D_{j}-\varepsilon) + D_{j}+\varepsilon \geq 0, $

\begin{equation*}
	|J_{j}(D_j, D_{j+1})| = \widetilde{V}_j(- \varepsilon, D_{j+1}- D_{j}-\varepsilon) + D_{j}+\varepsilon \leq |D_{j}| = O(1).
\end{equation*}

If $\; \widetilde{V}_j(-D_{j} - \varepsilon, D_{j+1}- 2 D_{j}-\varepsilon) + 2 D_{j}+\varepsilon < 0, $

\begin{equation*}
	|J_{j}(D_j, D_{j+1})| < \left|D_{j}\right|+ \dfrac{w_{j+1}}{w_{j}} \varepsilon = O(1).
\end{equation*}

\textbf{Case C. $D_j >0, D_{j+1}>0$}. In this case $T=\varepsilon > 0.$

\begin{equation*}
	J_j(D_j, D_{j+1})=\widetilde{V}_j(D_j+\varepsilon, D_{j+1}+\varepsilon) -\varepsilon,
\end{equation*}

We know that
\begin{align*}
	\widetilde{V}_j(D_{j} + \varepsilon, D_{j+1}+\varepsilon) &\geq D_{j} + \varepsilon,\\
	\widetilde{V}_j(D_{j} + \varepsilon, D_{j+1}+\varepsilon) &< \dfrac{1}{w_j}(D_{j}+\varepsilon).
\end{align*}

Thus, if 	$\; \widetilde{V}_j(D_{j}+\varepsilon, D_{j+1}+\varepsilon) -\varepsilon \geq 0, $

\begin{equation*}
	|J_{j}(D_j, D_{j+1})| < \dfrac{|D_{j}|}{w_j}+ \dfrac{w_{j+1}}{w_{j}} \varepsilon = O(1).
\end{equation*}

If  $\; \widetilde{V}_j(D_{j}+\varepsilon, D_{j+1}+\varepsilon) -\varepsilon < 0, $

\begin{equation*}
	|J_{j}(D_j, D_{j+1})| = \varepsilon - \widetilde{V}_j(D_{j}+\varepsilon, D_{j+1}+\varepsilon) \leq -D_{j} \:\text{which is not possible}.
\end{equation*}
\\
\textbf{Case D. $D_j < 0, \, D_{j+1}<0$}. In this case $T=-\varepsilon < 0.$

\begin{equation*}
	J_{j}(D_j, D_{j+1})=\widetilde{V}_j(D_{j}-\varepsilon, D_{j+1}-\varepsilon) +  \varepsilon,
\end{equation*}

We know that
\begin{align*}
	\widetilde{V}_j(D_{j}- \varepsilon, D_{j+1}-\varepsilon) &\leq D_{j}- \varepsilon,\\
	\widetilde{V}_j(D_{j}- \varepsilon, D_{j+1}-\varepsilon) &> \dfrac{1}{w_j}( D_{j}-\varepsilon).
\end{align*}

Thus, if $\; \widetilde{V}_j(D_{j}- \varepsilon, D_{j+1}- \varepsilon) + \varepsilon \geq 0, $

\begin{equation*}
	|J_{j}(D_j, D_{j+1})| = \widetilde{V}_j(D_{j}- \varepsilon, D_{j+1}-\varepsilon) + \varepsilon \leq D_{j} \:\text{which is not possible}.
\end{equation*}
\\

If $\; \widetilde{V}_j(D_{j} - \varepsilon, D_{j+1}- \varepsilon) + \varepsilon < 0, $

\begin{equation*}
	|J_{j}(D_j, D_{j+1})| < \dfrac{|D_{j}|}{w_j}+ \dfrac{w_{j+1}}{w_{j}} \varepsilon = O(1).
\end{equation*}

\section{Numerical experiment} \label{sec6}

\begin{figure}[!ht]
\begin{center}
\includegraphics[width=12cm]{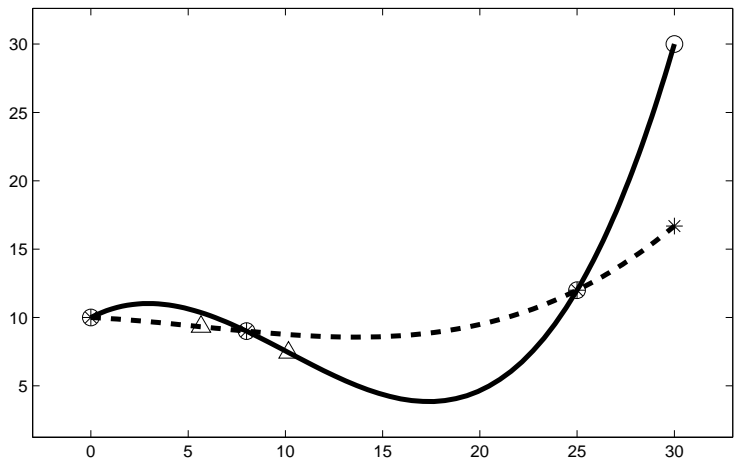}
\caption{\small\textit{In solid line: Lagrange polynomial, in dashed line: PPH polynomial.
Circles stand for Lagrange values at the nodes, asterisks for PPH and triangles for inflection points.}} \label{fig1}
\end{center}
\end{figure}

In this section we present three simple numerical experiments. The first one is dedicated to compare the convexity preservation between Lagrange and PPH reconstructions. Let us consider, the initial convex set of points, $(0,10)$, $(8,9)$, $(25,12)$ and $(30,30)\;$, that is $D_{j} >0, \, D_{j+1}>0 \;$. In Figure \ref{fig1} we have depicted the reconstruction operators corresponding to Lagrange and PPH and we have marked with triangles the inflection points for each reconstruction (5,66 and 10,16 respectively). We observe that PPH preserves convexity in a wider range $\frac{h_{j+1}+2 h_{j+2}}{6}=4.5,$ than Lagrange reconstruction does, see (\ref{eq:NUEconvexDif}).

Next experiment compares the approximation order of the considered reconstruction operators and allows us to check the effect on smooth functions of translation studied in section \ref{sec5}.
\\

Let be $X=(0,3,8,11,17,23,25,30,37,40)\dfrac{\pi}{20} \;$ a non uniform grid in $\, [0, 2\pi] \,$ and $f(x) = sin(x) \;$ a smooth function. Let us consider the set of initial points $\, (x_i, f(x_i)), \quad i=1,...,10. \,$ The approximation order $p \,$ obtained using Lagrange, PPH and PPH with translation, will be calculated for each reconstruction in a iterative way, dividing in two equal parts the intervals at each iteration $s,  \quad s= 0,1,...$.

The interval width at iteration $s \;$ will be $\approx \dfrac{h}{2^{s}}$, and the error\\
$E_{s}\approx C \left(\dfrac{h}{2^{s}}\right)^p$. So $\; \dfrac{E_{s-1}}{E_{s}} \approx 2^p \;$ and $p\approx \log_2\dfrac{E_{s-1}}{E_{s}}, \quad s= 1,2,...$.
\\

The error $E_{s}$ for each reconstruction $\,R(x) \,$ at iteration $\, s \,$ is calculated as a discrete approximation to $\; \left\|f(x)-R(x)\right\|_{\infty} \,$ evaluating  in a denser set of points.
\\

\begin{table}[!ht]
\begin{center}
\begin{tabular}{|c|c|c|c|c|}
\hline
 $ $  & $Lagrange$ & $PPH$ & $\dfrac{Translation}{\varepsilon=0.5}$ & $\dfrac{Translation}{\varepsilon=0.05}$\\[0.2cm]
\hline $s=1$ & $3.1461$ & $1.5701$ & $3.2622$ & $2.4126$ \\[0.2cm]
\hline $s=2$ & $3.7313$ & $2.9836$ & $3.5960$ & $3.3578$ \\[0.2cm]
\hline $s=3$ & $3.8978$ & $2.9959$ & $3.9280$ & $3.5412$ \\[0.2cm]
\hline $s=4$ & $3.9751$ & $2.9990$ & $3.9623$ & $3.7041$ \\[0.2cm]
\hline $s=5$ & $3.9938$ & $2.9997$ & $3.9811$ & $3.8264$ \\[0.2cm]
\hline
\end{tabular}
\caption{\small\textit{Approximation orders obtained at iteration $s,  \quad s= 0,1,..,5$.}} \label{tableorders}
\end{center}
\end{table}

\begin{figure}[!ht]
\begin{center}
\includegraphics[width=11cm]{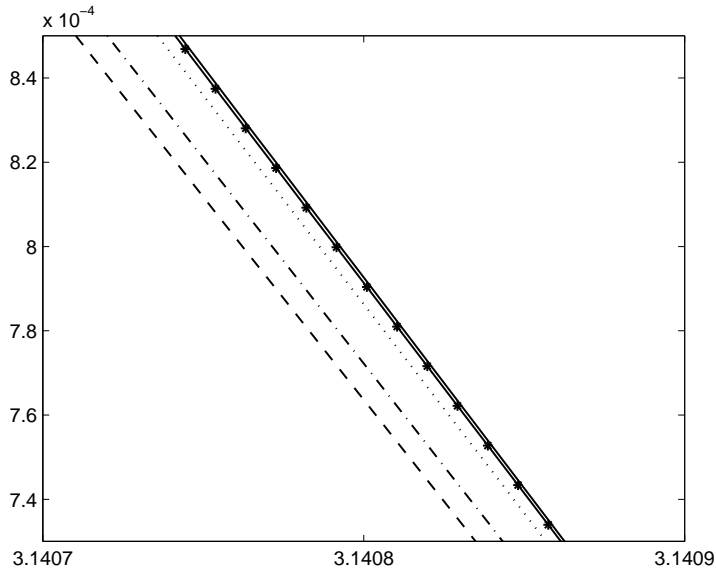}
\caption{\small\textit{In solid line: function $sin (x)$, in solid line with asterisks: Lagrange reconstruction, in dashed line: PPH reconstruction, in dash-dot line: Translated reconstruction for $\varepsilon= 0.05$, in dotted line: Translated reconstruction for $\varepsilon= 0.5$.}} \label{fig2}
\end{center}
\end{figure}

The approximation orders $p\,$ for each reconstruction are shown in Table \ref{tableorders}, where we can see, as it is well known, that Lagrange reconstruction attains fourth order accuracy with smooth functions. This same behavior would be desirable for the rest of the reconstructions in this case. However, PPH in this experiment obtains only third order accuracy and the reason is that the initial set of points contains inflection points of $\,sin (x)\,$ and near them, second order divided differences are of the same magnitude than $h$. Therefore, whether or not $D_{j}D_{j+1} \,$ is positive, from equation (\ref{eq:difMV}) it results that numerically $M_j-\widetilde{V}_j = O(h)$ and $f(x)-PPH_j(x) = O(h^3)$.
\\
This bad behavior of PPH in the vicinity of inflection points is corrected in the case of considering the translation version of PPH, as it can be seen in Table \ref{tableorders}.
\\
Because the proposed function is smooth, when $\,\varepsilon \,$ grows, translated reconstruction improves PPH order ($3.826$ for $\varepsilon = 0.05$, $3.981$ for $\varepsilon = 0.5$) and it get closer to the function as we can see in Figure \ref{fig2}, where we have depicted a zoom of the area around the inflection point $x= \pi$. We notice that outside from the inflection point areas, all four reconstructions are fourth order.
\\

Finally, a new experiment is proposed to analyze the adaptation to singularities of the considered reconstructions. The initial grid $\, X \,$ will be the same as in previous experiment and the proposed discontinuous function is
\begin{equation}
    g(x)=
        \left\{
            \begin{array}{lr}
                sin(x) & \text{if} \; x\leq 1.2 \pi,\\
                cos(x)+10 & otherwise,
            \end{array}
        \right\}
\end{equation}
thus, the set of initial points will be  $\, (x_i, g(x_i)), \quad i=1,...,10. \,$\\

\begin{figure}[!ht]
\begin{center}
\includegraphics[width=11cm]{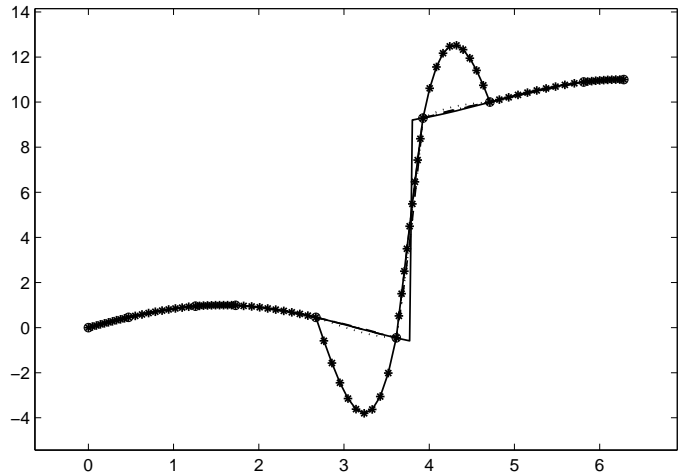}
\caption{\small\textit{In solid line: function $g(x)$, in solid line with asterisks Lagrange reconstruction, in dashed line: PPH reconstruction, in dash-dot line: Translated reconstruction for $\varepsilon= 0.05$, in dotted line: Translated reconstruction for $\varepsilon= 0.5$.}} \label{fig3}
\end{center}
\end{figure}

\begin{figure}[!ht]
\begin{center}
\includegraphics[width=11cm]{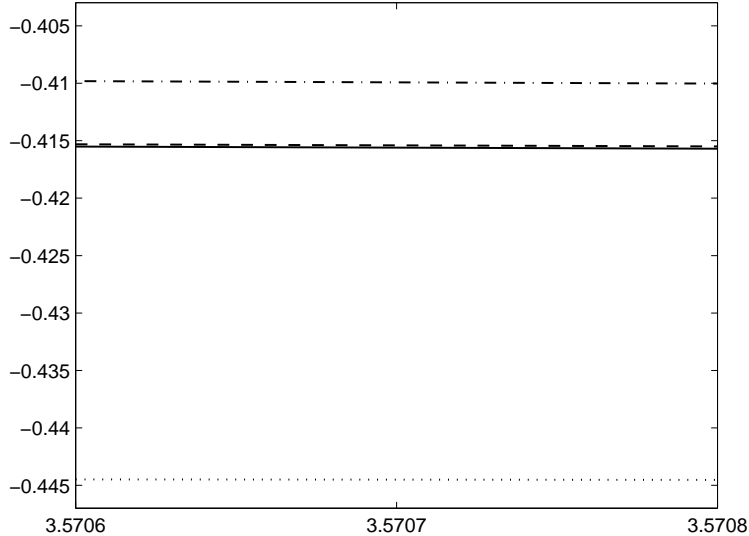}
\caption{\small\textit{In solid line: function $g(x)$, in dashed line: PPH reconstruction, in dash-dot line: Translated reconstruction for $\varepsilon= 0.05$, in dotted line: Translated reconstruction for $\varepsilon= 0.5$.}} \label{fig4}
\end{center}
\end{figure}

Figure \ref{fig3} shows the function $\, g(x) \,$ and the obtained reconstructions from the the initial points. We can see that around the singularity, Lagrange reconstruction looses the order and the Gibss phenomena appears. In this zone, PPH reconstruction is the best one and translations get closer and closer to PPH as $\varepsilon \,$ decreases, as it can also be seen more clearly in a zoom of this region in Figure \ref{fig4}.\\
Equation (\ref{eq:bNUETransP2}) allows to ensure that translation get close to PPH with decreasing $\, \varepsilon. \,$ However it might happen that for some larger $\, \varepsilon \,$ and for some regions the translation reconstruction get close to PPH that for some others smaller values of $\, \varepsilon, \,$ what does not enter in contradiction with (\ref{eq:bNUETransP2}).
\\

\section{Conclusions}\label{sec7}

We have studied the PPH reconstruction operator over non uniform grids. For this purpose, the arithmetic and harmonic means used in the uniform case have been changed for weighted means. We have seen that coefficients for reconstruction operators over non uniform grids are similar to those in the uniform case, changing the original mean by its weighted counterpart.

We have also studied some properties of the new reconstruction operator, checking that its behavior improves comparing with Lagrange near a singularity. For smooth convex functions, the PPH reconstruction gets fourth order accuracy as it must be. In the vicinity of an inflection point the order reduces. To avoid this fact, we have proposed a new mean making use of a translation strategy. This way, we maintain fourth order accuracy with smooth functions, convex or not, at the same time that we offer adaptation near singularities.

Also we have obtained the conditions under which Lagrange and PPH reconstruction operators preserve convexity.

Finally we have carried out some simple numerical experiments to reinforce the theoretical results.

\clearpage

\end{document}